# THE BERNOULLI SIEVE REVISITED

By Alexander V. Gnedin, Alexander M. Iksanov,[1]
Pavlo Negadajlov[1] and Uwe Rösler[1]

*Utrecht University, University of Kiev, University of Kiev and University of Kiel*

We consider an occupancy scheme in which "balls" are identified with $n$ points sampled from the standard exponential distribution, while the role of "boxes" is played by the spacings induced by an independent random walk with positive and nonlattice steps. We discuss the asymptotic behavior of five quantities: the index $K_n^*$ of the last occupied box, the number $K_n$ of occupied boxes, the number $K_{n,0}$ of empty boxes whose index is at most $K_n^*$, the index $W_n$ of the first empty box and the number of balls $Z_n$ in the last occupied box. It is shown that the limiting distribution of properly scaled and centered $K_n^*$ coincides with that of the number of renewals not exceeding $\log n$. A similar result is shown for $K_n$ and $W_n$ under a side condition that prevents occurrence of very small boxes. The condition also ensures that $K_{n,0}$ converges in distribution. Limiting results for $Z_n$ are established under an assumption of regular variation.

**1. Introduction.** The Bernoulli sieve is a simple recursive allocation of $n$ "balls" in infinitely many "boxes" indexed $1, 2, \ldots$. Let $\xi_1, \xi_2, \ldots$ be random values sampled independently from a given probability distribution on $(0, 1)$. At the first step each of $n$ balls is dropped in box 1 with probability $\xi_1$. At the second step each of the remaining balls is dropped in box 2 with probability $\xi_2$, and so on. The procedure is iterated until all $n$ balls get allocated. It is easy to see that the probability that a particular ball lands in box $j$ is equal to

$$P_j = \bar{\xi}_1 \cdots \bar{\xi}_{j-1} \xi_j, \qquad i \in \mathbb{N}, \tag{1}$$

here and hereafter $\bar{x} := 1 - x$.

Received April 2008.
[1]Supported by the German Scientific Foundation project no. 436UKR 113/93/0-1.

Random discrete probability distributions with frequencies $(P_j)$ of the form (1) are called *residual allocation* or *stick-breaking* models [1, 4, 24]. For instance, in the most popular and analytically best tractable case $(P_j)$ follows the GEM (Griffiths–Engen–McCloskey) distribution, which appears when the law of $\xi_1$ is $\mathrm{beta}(1,\theta)$ with some $\theta > 0$.

Note that $(P_j)$ is a (nonrandom) geometric distribution if the law of $\xi_1$ is a Dirac mass $\delta_x$ at some point $x \in (0,1)$; we shall exclude this and some other cases by the assumption that the support of $\xi_1$ is not a set like $\{1 - x^j : j \in \mathbb{N}_0 := \mathbb{N} \cup \{0\}\}$. See [12] for a survey of results on sampling from nonrandom discrete distributions with infinitely many positive masses.

A random combinatorial structure which captures the occupancy of boxes is the *weak composition* $C_n^*$ comprised of nonnegative integer *parts* summing up to $n$. We speak of *weak* composition meaning that zero parts are allowed, for instance, the sequence $(2,3,0,1,0,0,1,0,0,0,\ldots)$ (padded by infinitely many 0's) is a possible value of $C_7^*$. Two related structures which contain less information are *composition* $C_n$ of $n$ obtained by discarding zero parts of $C_n^*$ and a *partition* of $n$ obtained by arranging the parts in nonincreasing order [these are $(2,3,1,1)$ and $(3,2,1,1)$, respectively, in the example]. In the GEM case, the law of the partition is widely known as Ewens' sampling formula (ESF), and the law of composition is a size-biased version of the ESF; see [1, 24].

Functionals of $C_n^*$ studied in this paper are as follows:

$K_n$ the number of boxes occupied by at least one ball,
$K_n^*$ the largest index of occupied box,
$K_{n,0} = K_n^* - K_n$ the number of empty boxes with index not exceeding $K_n^*$,
$W_n$ the index of the first empty box,
$Z_n$ the number of balls in the last box.

For $r = 1, 2, \ldots, n$, we also denote by $K_{n,r}$ the number of boxes occupied by exactly $r$ balls.

In [10] it was observed that $K_n$ can be studied by tools of the renewal theory, and it was shown that under certain moment conditions the distribution of $K_n$ is asymptotically normal. The composition $C_n$ in this case has some common features with logarithmic combinatorial structures [1]; in particular, $K_n$ exhibits a logarithmic growth.

Throughout, we shall also rely on the following alternative construction of $C_n^*$'s. Let $\{S_k : k \in \mathbb{N}_0\}$ be a zero-delayed random walk with a step distributed like $(-\log \bar{\xi}_1)$. For $E_1, E_2, \ldots$ an independent random sample from the standard exponential distribution, also independent of $(S_k)$, think of the event $E_j \in (S_{k-1}, S_k)$ as a ball dropped in box $k$. A composition of $n$ is defined as the sequence of occupancy numbers in the natural order of intervals $(S_{k-1}, S_k), k = 1, 2, \ldots$. In what follows we will often use $E_{1,n} \leq E_{2,n} \leq \cdots \leq E_{n,n}$ the order statistics of $E_1, \ldots, E_n$.

The equivalence with the Bernoulli sieve construction is established via the mapping $y \mapsto e^{-y}, y > 0$, which allows to identify $\bar{\xi}_1 \cdots \bar{\xi}_k$ (we tacitly assume that this equals 1 for $k = 0$) with $\exp(-S_k), k \in \mathbb{N}_0$, transforms $(S_{k-1}, S_k)$ into interval of size $P_k$ and transforms $(E_1, \ldots, E_n)$ in a uniform $[0,1]$ sample. By this transformation, the event $E_j \in (S_{k-1}, S_k)$ occurs when the $j$th coordinate of the uniform sample falls in the $k$th interval $(\exp(-S_k), \exp(-S_{k-1}))$. This works, because a point sampled from the uniform $[0,1]$ distribution falls in the $k$th interval with probability $P_k$.

As $n$ varies, compositions $C_n^*$ satisfy the following consistency conditions:

- (SC) Sampling consistency: if one of $n$ balls is chosen uniformly at random and removed from the box it occupies, the resulting weak composition has the same probability law as $C_{n-1}^*$.
- (DP) Deletion property: if the first box is inspected and it turns out that it contains $k$ balls, then deleting the first box[2] yields a weak composition with the same probability law as $C_{n-k}^*$.

Condition (SC) follows from the independence of $(S_k)$ and $(E_j)$ and exchangeability of $(E_j)$, and condition (DP) follows from the renewal property of $(S_k)$ and the memoryless property of the exponential distribution. Both conditions also hold for the associated compositions, which means that the sequence $(C_n)$ is a *regenerative composition structure*, as introduced in [14]. Note that $K_n, K_{n,r}$ are, in fact, functionals of the *partition structure* [4, 24] which is obtained by discarding the ordering of parts in the $C_n$'s.

The random walk $(S_k)$ can be viewed as the range of a compound Poisson process, that is, a subordinator $\{T_s : s \geq 0\}$ whose Lévy measure is the distribution of $(-\log \bar{\xi})$. One obtains a larger class of composition structures by considering a general zero-drift subordinator and using the open gaps comprising the complement of its range in the role of boxes. It is known that normal limits for the number of parts $K_n$ are typical for regenerative composition structures whose Lévy measure is infinite and has the right tail slowly varying at 0 [3, 17], although $K_n$ exhibits then growth faster than logarithmic. The limits are no longer normal if the right tail of Lévy measure is regularly varying at 0 with positive index, as, for example, it is the case for stable subordinators [16].[3]

In this paper we dwell on the case of the Bernoulli sieve and obtain considerable extensions of the results of [10]. In particular, we derive an exhaustive criterion for the existence of limiting distribution of properly normalized and centered $K_n^*$ (Theorem 2.1). Then, under a side condition, we do the same for $K_n$ and $W_n$ (Theorem 2.3). Among other things, this condition ensures our most delicate result, which states that $K_{n,0}$ converges in distribution directly, without centering or scaling (Theorem 2.2). A similar result also holds for $K_{n,0} + K_{n,1}$ (Proposition 5.2). In the GEM case, the limiting law of $K_{n,0}$ is mixed Poisson (Proposition 5.1). Asymptotic properties of $Z_n$ are revealed in Theorem 2.4.

---

[2] In the example, the elimination transforms $(2, 3, 0, 1, 0, 0, 1, 0, 0, 0, \ldots)$ to $(3, 0, 1, 0, 0, 1, 0, 0, 0, \ldots)$.

[3] It should be noticed that in the case of infinite Lévy measure the closed range of subordinator is a Cantor set, thus, with positive probability the set of empty boxes between $E_{j-1,n}$ and $E_{j,n}$ is infinite.

The rest of the paper is organized as follows. In Section 2 we formulate our principal results and give examples. In Section 3 we extend the idea of representing a regenerative composition via a Markov chain [14, 15] to the weak compositions. We also collect necessary distributional recursions. Other sections are devoted to study of particular functionals. The Appendix summarizes some known asymptotic results about the number of renewals.

**2. Main results.** Consider the process that counts renewals

$$(2) \qquad N_t := \inf\{k \geq 1 : S_k \geq t\}, \qquad t \geq 0.$$

Our idea is to connect possible convergence in distribution of $(N_{\log n} - b_n)/a_n$ to some nondegenerate and proper probability law with the convergence of $(K_n^* - b_n)/a_n$, $(K_n - b_n)/a_n$, $(K_n - K_{n,1} - b_n)/a_n$ and $(W_n - b_n)/a_n$ to the same law. The first connection can be anticipated in the view of identity $K_n^* = N_{E_{n,n}}$ and by recalling the fact from the extreme-value theory that $E_{n,n} - \log n$ has a limiting distribution (of Gumbel type).

Introduce the moments

$$\mu := \mathbb{E}(-\log \bar{\xi}) \quad \text{and} \quad \sigma^2 := \mathrm{var}(-\log \bar{\xi}),$$

which may be finite or infinite.

THEOREM 2.1.  *The following assertions are equivalent:*

(i) *There exist constants $\{a_n, b_n : n \in \mathbb{N}\}$ with $a_n > 0$ and $b_n \in \mathbb{R}$ such that, as $n \to \infty$, the variable $(K_n^* - b_n)/a_n$ converges weakly to some nondegenerate and proper distribution.*

(ii) *The distribution of $(-\log \bar{\xi})$ either belongs to the domain of attraction of a stable law, or the function $\mathbb{P}\{-\log \bar{\xi} > x\}$ slowly varies at $\infty$.*

*Furthermore, this limiting distribution is as follows:*

(a) *If $\sigma^2 < \infty$, then for $b_n = \mu^{-1} \log n$ and $a_n = (\mu^{-3} \sigma^2 \log n)^{1/2}$, the limiting distribution is standard normal.*
(b) *If $\sigma^2 = \infty$ and*

$$\int_x^1 (\log y)^2 \mathbb{P}\{\bar{\xi} \in dy\} \sim L(-\log x) \qquad \text{as } x \to 0,$$

*for some $L$ slowly varying at $\infty$, then for*

$$b_n = \mu^{-1} \log n, a_n = \mu^{-3/2} c_{[\log n]}$$

*and $c_n$ any sequence satisfying $\lim_{n\to\infty} nL(c_n)/c_n^2 = 1$, the limiting distribution is standard normal.*

(c) *Assume that the relation*

(3)  $$\mathbb{P}\{\bar{\xi} \leq x\} \sim (-\log x)^{-\alpha} L(-\log x) \qquad \text{as } x \to 0$$

*holds with $L$ slowly varying at $\infty$ and $\alpha \in [1,2)$, and assume that $\mu < \infty$ if $\alpha = 1$, then for $b_n = \mu^{-1} \log n$, $a_n = \mu^{-(\alpha+1)/\alpha} c_{[\log n]}$ and $c_n$ any sequence satisfying $\lim_{n\to\infty} nL(c_n)/c_n^\alpha = 1$, the limiting distribution is $\alpha$-stable with characteristic function*

$$t \mapsto \exp\{-|t|^\alpha \Gamma(1-\alpha)(\cos(\pi\alpha/2) + i\sin(\pi\alpha/2)\operatorname{sgn}(t))\}, \qquad t \in \mathbb{R}.$$

(d) *Assume that $\mu = \infty$ and the relation (3) holds with $\alpha = 1$. Let $c$ be any positive function satisfying $\lim_{x\to\infty} xL(c(x))/c(x) = 1$ and set $\psi(x) := x \int_{\exp(-c(x))}^1 \mathbb{P}\{\bar{\xi} \leq y\}/y\, dy$. Let $b$ be any positive function satisfying $b(\psi(x)) \sim \psi(b(x)) \sim x$. Then, with $b_n = b(\log n)$ and $a_n = b(\log n) \times c(b(\log n))/\log n$, the limiting distribution is $1$-stable with characteristic function*

(4) $$t \mapsto \exp\{-|t|(\pi/2 - i\log|t|\operatorname{sgn}(t))\}, \qquad t \in \mathbb{R}.$$

(e) *If the relation (3) holds with $\alpha \in [0,1)$, then, for $b_n = 0$ and $a_n := \log^\alpha n/L(\log n)$, the limiting distribution is the scaled Mittag–Leffler law $\theta_\alpha$ (exponential, if $\alpha = 0$) characterized by the moments*

$$\int_0^\infty x^n \theta_\alpha(dx) = \frac{n!}{\Gamma^n(1-\alpha)\Gamma(1+n\alpha)}, \qquad n \in \mathbb{N}.$$

Asymptotic analysis of the number of empty boxes $K_{n,0}$ involves

$$\nu := \mathbb{E}(-\log \xi).$$

Our next result determines explicitly the limiting distribution of $K_{n,0}$.

THEOREM 2.2. *For $n \to \infty$, $K_{n,0}$ has the following asymptotic properties:*

(a) *If $\nu < \infty$, then $K_{n,0}$ converges in distribution to a random variable $K_{\infty,0}$. If also $\mu < \infty$, then*

$$\mathbb{P}\{K_{\infty,0} \geq i\} = \frac{1}{\mu} \sum_{j=1}^\infty \frac{\mathbb{E}\bar{\xi}^j}{j} \mathbb{P}\{K_{j,0} = i-1\}, \qquad i \in \mathbb{N},$$

*and $\mathbb{E}K_{\infty,0} = \nu/\mu$, but if $\mu = \infty$, then $K_{\infty,0} = 0$ a.s.*

(b) *Assume that for some $\delta > 0$ both*

(5) $$\mathbb{E}\bar{\xi}^{-\delta} < \infty \quad \text{and} \quad \mathbb{E}\xi^{-\delta} < \infty,$$

then

$$\lim_{n \to \infty} \mathbb{E}K_{n,0} = \nu/\mu \in (0, \infty).$$

*On the other hand, if $\nu = \infty$ and $\mu < \infty$, then $\lim_{n \to \infty} \mathbb{E}K_{n,0} = \infty$.*

If $\mu < \infty$, the limiting variable $K_{0,\infty}$ has interpretation in terms of a model with infinitely many "balls" and "boxes" [13]. Specifically, one can take gaps between consecutive points in a stationary renewal process on $\mathbb{R}$ in the role of "boxes," and points of an independent Poisson process with the intensity measure $e^{-x} dx$ ($x \in \mathbb{R}$) in the role of "balls." Other functionals of $C_n^*$ like $K_{n,r}$ the number of parts equal $r$ also have limiting forms realizable in the infinite model.

By virtue of $K_n = K_n^* - K_{n,0}$ and because Theorem 2.1 implies that $a_n \to \infty$, one can conclude that boundedness of $K_{n,0}$ ($K_{n,1}$) in probability would lead to the following fact: if $(K_n^* - b_n)/a_n$ weakly converges to some proper probability law, then $(K_n - b_n)/a_n$ (($K_n - K_{n,1} - b_n)/a_n$) weakly converges to the same law. According to Theorem 2.2 (Proposition 5.2), the condition $\nu < \infty$ ensures even a more delicate property that $K_{n,0}$ ($K_{n,0} + K_{n,1}$) converges in distribution. Similar argument applies to $W_n$ and leads to the next result.

THEOREM 2.3.  *If $\nu < \infty$, then all the assertions of Theorem 2.1 remain valid with $K_n^*$ replaced by $K_n$, $K_n - K_{n,1}$ or $W_n$.*

Under the condition $\sigma^2 < \infty$, the normal limit for $K_n$ was established in [10], Proposition 10, by a method which required asymptotic expansion of moments. A generalization for a larger class of random compositions appeared in [11], Proposition 8.

THEOREM 2.4.  *If $\mu < \infty$, then $Z_n \overset{d}{\to} Z$ as $n \to \infty$, where $Z$ has distribution*

$$\mathbb{P}\{Z = k\} = \frac{\mathbb{E}\xi^k}{\mu k}, \qquad k = 1, 2, \ldots.$$

*If (3) holds with some $\alpha \in [0,1)$, then for $\alpha \in (0,1)$*

$$\frac{\log Z_n}{\log n} \overset{d}{\to} \text{beta}(1-\alpha, \alpha)$$

and the degenerate limit distribution $\delta_1$ appears for $\alpha = 0$. If (3) holds with $\alpha = 1$ and if $\mu = \infty$, then with $m(x) := \int_0^x \mathbb{P}\{-\log \bar\xi > y\}\,dy$ we have the convergence

$$\frac{m(\log Z_n)}{m(\log n)} \xrightarrow{d} \mathrm{uniform}[0,1].$$

In the examples to follow $X_n$ stands for any of the variables $K_n^*$, $K_n$, $K_n - K_{n,1}$ or $W_n$.

EXAMPLE 2.5. Assume $\bar\xi$ has a beta$(c,b)$ density

$$\mathbb{P}\{\bar\xi \in dx\} = \frac{x^{c-1}(1-x)^{b-1}}{B(c,b)}\,dx, \qquad x \in [0,1],$$

with some $c,b > 0$ and $B(\cdot,\cdot)$ denoting the beta function [hence, the law of $\xi_1$ is beta$(b,c)$]. In this case the moments are finite and given by

$$\mu = \Psi(c+b) - \Psi(c),$$
$$\nu = \Psi(c+b) - \Psi(b),$$
$$\sigma^2 = \Psi'(c) - \Psi'(c+b),$$

where $\Psi(x) = \Gamma'(x)/\Gamma(x)$ denotes the logarithmic derivative of the gamma function. Therefore, as $n \to \infty$,

$$\frac{X_n - \mu^{-1}\log n}{(\mu^{-3}\sigma^2 \log n)^{1/2}} \xrightarrow{d} \mathrm{normal}(0,1).$$

Above that, $Z_n \xrightarrow{d} Z$ with $Z$ having distribution

$$\mathbb{P}\{Z = k\} = \frac{\Gamma(c+b)}{\mu \Gamma(b)}\frac{\Gamma(k+b)}{k\Gamma(k+b+c)}, \qquad k = 1,2,\ldots.$$

The number of empty boxes $K_{n,0}$ converges in distribution and in the mean to a random variable $K_{\infty,0}$ with some nondegenerate distribution. For $b \neq 1$ an explicit form of the limiting distribution is still a challenge. For $b = 1$ Proposition 5.1 gives the generating function

$$\mathbb{E} s^{K_{\infty,0}} = \frac{\Gamma(1+c)\Gamma(1+c-cs)}{\Gamma(1+2c-cs)}, \qquad s \in [0,1].$$

In particular, for integer $c$ the distribution of $K_{\infty,0}$ is the convolution of $c$ geometric distributions with parameters $k^{-1}(k+c)$, $k = 1,2,\ldots,c$.

EXAMPLE 2.6. Suppose $\bar\xi$ has distribution function

$$\mathbb{P}\{\bar\xi \leq x\} = \frac{1}{1 - \log x}, \qquad x \in (0,1).$$

Then the condition (3) holds with $\alpha = 1$ and $\mu = \infty$. Since

$$\mathbb{P}\{-\log \xi > x\} = \frac{-\log(1 - e^{-x})}{1 - \log(1 - e^{-x})}, \qquad x \geq 0,$$

for $x \to \infty$ we have $\mathbb{P}\{-\log \xi > x\} \sim e^{-x}$, whence $\nu < \infty$. Therefore, as $n \to \infty$,

(6) $$\frac{(\log \log n)^2}{\log n} X_n - \log \log n - \log \log \log n$$

converges in distribution to the spectrally negative 1-stable law with characteristic function (4). Since $\mathbb{P}\{-\log \bar{\xi} > x\} = (x+1)^{-1}$ holds for $x > 0$, the normalizing constants in (6) can be calculated in the same way as in [20], Proposition 2. Above that,

$$\frac{\log \log Z_n}{\log \log n} \xrightarrow{d} \text{uniform}[0,1] \quad \text{and} \quad K_{n,0} \xrightarrow{d} \delta_0.$$

EXAMPLE 2.7. For $\bar{\xi}$ with distribution

$$\mathbb{P}\{\bar{\xi} \leq x\} = \frac{-\log(1-x)}{1 - \log(1-x)}, \qquad x \in (0,1),$$

we have $\sigma^2 < \infty$ but $\nu = \infty$, hence, Theorem 2.3 is not applicable.

**3. Compositions, Markov chains and recursions.** Weak composition $C_n^*$ can be identified with a path of a time-homogeneous nonincreasing Markov chain $Q_n^*$ on integers which start at $n$, terminate at 0 and have nonnegative integer decrements equal to the parts of $C_n^*$. Similarly, composition $C_n$ can be identified with the path of a Markov chain $Q_n$, whose decrements are positive until absorbtion at state 0. For fixed $n$, in terms of "balls-in-boxes," $Q_n^*(k)$ is the number of exponential points which fall outside the first $k$ spacings induced by $(S_i)$, and $Q_n(k)$ is the number of exponential points which fall outside the first $k$ spacings containing at least one of $E_j$, for $k = 0, 1, \ldots$.

Following terminology from [14], the transition probabilities are determined by the *decrement matrices*

(7) $$q^*(n:m) := \binom{n}{m} \mathbb{E}[\xi^{n-m} \bar{\xi}^m], \qquad m = 0, \ldots, n,$$

(8) $$q(n:m) := \binom{n}{m} \frac{\mathbb{E}[\xi^{n-m} \bar{\xi}^m]}{1 - \mathbb{E}[\bar{\xi}^n]}, \qquad m = 1, \ldots, n,$$

which specify the probability distribution of the first part of $C_n^*$, respectively, $C_n$. By this representation, $q^*(n:m)$ and $q(n:m)$ are the transition probabilities from $n \geq 0$ to $n-m$ for the Markov chains $Q_n^*$ and $Q_n$, respectively.

Introduce the total frequency of boxes whose indices are larger than $j$,

$$\zeta_j := \bar{\xi}_1 \cdots \bar{\xi}_j.$$

From the construction of $C_n^*$ it is clear that

(9) $$\mathbb{P}\{Q_n^*(j) = n - m\} = \binom{n}{m} \mathbb{E}[\zeta_j^{n-m}(1-\zeta_j)^m],$$

which is the multistep generalization of (7). Also,

(10) $$\mathbb{P}\{K_n^* > k\} = \mathbb{P}\{Q_n^*(k) > 0\} = \mathbb{E}[1 - (1-\zeta_k)^n].$$

The variables we are interested in have obvious interpretation via the Markov chains. Thus, the absorbtion time $Q_n^*$, that is, the number of steps the chain needs to approach 0 is $K_n^*$, and the absorption time of $Q_n$ is $K_n$.

The Markov property leads to distributional recursions

(11) $$K_0^* := 0, \qquad K_n^* \stackrel{d}{=} K_{A_n^*}^* + 1, \qquad n \in \mathbb{N},$$

and

(12) $$K_0 = 0, \qquad K_n \stackrel{d}{=} K_{A_n} + 1, \qquad n \in \mathbb{N},$$

where $A_n^*$ is assumed independent of $\{K_j^* : j \in \mathbb{N}\}$ and $A_n^* \stackrel{d}{=} Q_n^*(1)$; $A_n$ is independent of $\{K_j : j \in \mathbb{N}\}$ and $A_n \stackrel{d}{=} Q_n(1)$. So the law of $A_n^*$ is $q^*(n:\cdot)$ and the law of $A_n$ is $q(n:\cdot)$.

Now let $V_n$ be the number of balls that fall to the right from the first empty box. For instance, for weak compositions $(1, 2, 1, 0, 2, 0, 0, 3, 0, 0, \ldots)$, $(1, 2, 1, 2, 3, 0, 0, \ldots)$ and $(0, 1, 2, 1, 2, 3, 0, 0, \ldots)$, the value of $V_9$ is $5, 0$ and $9$, respectively. Then

(13) $$K_{0,0} = 0, \qquad K_{n,0} \stackrel{d}{=} K_{V_n,0} + 1_{\{V_n > 0\}}, \qquad n \in \mathbb{N},$$

where on the right-hand side $V_n$ is independent of $\{K_{n,0} : n \in \mathbb{N}\}$. Here and below, $1_{\{\cdots\}}$ is 1 if $\cdots$ holds true and is 0 otherwise. Furthermore,

(14) $$K_0^* = 0, \qquad K_n^* \stackrel{d}{=} K_{V_n}^* + W_n - 1_{\{V_n = 0\}}, \qquad n \in \mathbb{N},$$

where on the right-hand side $(V_n, W_n)$ are independent of $\{K_n^* : n \in \mathbb{N}\}$. Finally, we remark that

(15) $$K_{n,0} \stackrel{d}{=} K_{A_n^*, 0} + 1_{\{A_n^* = n\}}, \qquad n \in \mathbb{N},$$

where $A_n^*$ is independent of $\{K_{n,0} : n \in \mathbb{N}\}$.

The chain $Q_n$ visits a given state $m$ with the same probability as $Q_n^*$. We denote this probability by

$$g(n,m) := \mathbb{P}\{Q_n^*(0) = m\} + \sum_{j=1}^{\infty} \mathbb{P}\{Q_n^*(j) = m, Q_n^*(j-1) \neq m\}$$

$$= \sum_{j=0}^{\infty} \mathbb{P}\{Q_n(j) = m\}.$$

In principle, $g$ can be computed from (9), but there is a simpler and more general formula which involves only $\mathbb{E}[1-\bar{\xi}^k], k = m, \ldots, n$; see [14], Theorem 9.2.

It was shown in [10], Proposition 5, that under the assumption $\mu < \infty$

(16) $$\lim_{n \to \infty} g(n,m) = \frac{1 - \mathbb{E}\bar{\xi}^m}{\mu m}$$

and the same argument[4] allows one to show that $\lim_{n \to \infty} g(n,m) = 0$ if $\mu = \infty$.

**4. Index of the last occupied box and a proof of Theorem 2.1.** By results of the renewal theory (which we summarize in Proposition A.1), it is enough to show the equivalence

(17) $$\frac{K_n^* - b_n}{a_n} \xrightarrow{d} X \iff \frac{N_{\log n} - b_n}{a_n} \xrightarrow{d} X,$$

where $X$ is a random variable with a proper and nondegenerate probability distribution $F$, and $N_t$ is as in (2). Assuming that the convergence in the left-hand side of (17) holds with $a_n \to \infty$, we have for $y > 0$

$$\mathbb{P}\left\{\frac{K_n^* - b_n}{a_n} > x\right\} = \mathbb{P}\left\{\frac{N_{E_{n,n}} - b_n}{a_n} > x\right\}$$

$$\leq \mathbb{P}\left\{\frac{N_{\log n + y} - b_n}{a_n} > x\right\}\mathbb{P}\{E_{n,n} - \log n \leq y\}$$

$$+ \mathbb{P}\{E_{n,n} - \log n > y\}.$$

By subadditivity of the number of renewals, $N_{\log n + y}$ does not exceed stochastically the sum $N_{\log n} + N_y'$ with independent terms and $N_y' \stackrel{d}{=} N_y$, hence, we can estimate the above further as

$$\leq \mathbb{P}\left\{\frac{N_{\log n} - b_n}{a_n} + \frac{N_y'}{a_n} > x\right\}\mathbb{P}\{E_{n,n} - \log n \leq y\} + \mathbb{P}\{E_{n,n} - \log n > y\}.$$

---
[4]The proof on top of [10], page 86, must be corrected by changing $n-m$ to $n-m+1$.

By the selection principle, there exists an increasing subsequence $(n_k)$ such that the variable $(N_{\log n_k} - b_{n_k})/a_{n_k}$ converges weakly to some measure $F'$, say. Recalling the convergence of $E_{n,n} - \log n$ and sending $y$ to $\infty$, we have $F(x, \infty) \leq F'(x, \infty)$ at all joint continuity points of $F(x, \infty)$ and $F'(x, \infty)$.

Similarly, for $y < 0$,

$$\mathbb{P}\left\{\frac{N_{E_{n,n}} - b_n}{a_n} > x\right\} \geq \mathbb{P}\left\{\frac{N_{\log n + y} - b_n}{a_n} > x\right\}\mathbb{P}\{E_{n,n} - \log n > y\}$$

$$\geq \mathbb{P}\left\{\frac{N_{\log n} - b_n}{a_n} - \frac{N'_{-y}}{a_n} > x\right\}\mathbb{P}\{E_{n,n} - \log n > y\}.$$

Letting again $n \to \infty$ along $(n_k)$ and then sending $y$ to $-\infty$, we have $F(x, \infty) \geq F'(x, \infty)$ at all joint continuity points of $F$ and $F'$. Therefore, $F = F'$ and since the limit does not depend on subsequence, we conclude that $(N_{\log n} - b_n)/a_n \overset{d}{\to} X$.

Obviously the number of balls outside the first box, $A_n^*$, goes to $\infty$ a.s., which together with an application of (11) implies that $K_n^*$ cannot converge in distribution if no scaling or centering is imposed. Nor can $K_n^* - b_n$, for any unbounded sequence $b_n > 0$. Indeed, if the convergence were the case, from the convergence of $E_{n,n} - \log n$ and a.s. monotonicity of $N_t$ would follow that $N_{\log n} - b_n$ were bounded in probability, which is known to be false. Following the same line of argument, one can prove that $(K_n^* - b_n)/a_n$ also cannot converge in distribution if $a_n$ is either bounded or unbounded but does not go to $\infty$.

To establish the result in the reverse direction, we prefer to exploit the multiplicative form of renewal process. For each $\varepsilon > 0$ define

$$M_n^{(\varepsilon)} := \inf\{k \geq 1 : n\bar{\xi}_1\bar{\xi}_2\cdots\bar{\xi}_k \leq \varepsilon\}, \qquad n \in \mathbb{N},$$

and notice that $M_n^{(1)} = N_{\log n}$. Assume that $(M_n^{(1)} - b_n)/a_n \overset{d}{\to} X$, where $X$ is a random variable with a proper and nondegenerate distribution $F$. By Proposition A.1 from the Appendix, $F$ is continuous and $a_n$ slowly varies. Also, Proposition A.1 provides an explicit form of $b_n$. Using this, we conclude that $(M_n^{(\varepsilon)} - b_n)/a_n \overset{d}{\to} X$, no matter what $\varepsilon$ is.

For fixed $x \in \mathbb{R}$ and $n$ sufficiently large put $k_n := \lfloor a_n x + b_n \rfloor$. Since for large $n$

$$\mathbb{E}[1 - (1 - \zeta_{k_n})^n] \geq \mathbb{E}[1 - (1 - \zeta_{k_n})^n \mathbf{1}(\zeta_{k_n} > \varepsilon/n)]$$

$$\geq (1 - (1 - \varepsilon/n)^n)\mathbb{P}\{M_n^{(\varepsilon)} > k_n\},$$

letting in (10) first $n \to \infty$ and then $\varepsilon \to 0$, we obtain $\liminf_{n \to \infty} \mathbb{P}\{K_n^* > k_n\} \geq F(x, \infty)$. On the other hand, for large $n$,

$$\mathbb{E}[1 - (1 - \zeta_{k_n})^n] \leq (1 - (1 - \varepsilon/n)^n)\mathbb{P}\{M_n^{(\varepsilon)} \leq k_n\} + \mathbb{P}\{M_n^{(\varepsilon)} > k_n\}.$$

Sending in (10) first $n \to \infty$ and then $\varepsilon \to \infty$, we obtain $\limsup_{n\to\infty} \mathbb{P}\{K_n^* > k_n\} \leq F(x, \infty)$. Combining the lower and upper limit, we conclude that $(K_n^* - b_n)/a_n \xrightarrow{d} X$, as desired.

**5. The number of empty boxes and a proof of Theorem 2.2.** In the setting of GEM distribution, that is, when $\xi \stackrel{d}{=} \text{beta}(1, \theta)$, the distribution of the number $K_{n,r}$ of boxes occupied by exactly $r$ balls is asymptotically Poisson$(\theta/r)$, for every $r > 0$. See, for example, [1], Theorem 4.17, where the fact appears in connection with the cycle structure of random $\theta$-biased permutations. Quite unexpectedly, the limit law of $K_{n,0}$ is not Poisson. In the spirit of discussion after Theorem 2.2, the limit variable may be interpreted in terms of the Poisson process $\Pi_1$ (boxes) of intensity $(\theta/x)\, dx, x \in \mathbb{R}_+$, and another independent rate-1 Poisson process $\Pi_2$ (balls) on $\mathbb{R}_+$: $K_{0,\infty}$ is the number of gaps in $\Pi_1$ that are to the right of the leftmost atom of $\Pi_2$ and do not contain points of $\Pi_2$.

PROPOSITION 5.1. *If $\xi$ has* beta$(1, \theta)$ *distribution* $(\theta > 0)$, *then $K_{n,0}$ converges in distribution to a variable $K_{\infty,0}$ with*

$$(18) \qquad \mathbb{E}s^{K_{\infty,0}} = \frac{\Gamma(1+\theta)\Gamma(1+\theta-\theta s)}{\Gamma(1+2\theta-\theta s)}, \qquad s \in [0,1],$$

*which is the generating function of a mixed Poisson distribution with random parameter $\theta|\log \xi|$.*

PROOF. For $j = 1, 2, \ldots$ let $M_j \stackrel{d}{=} \text{geometric}(j/(\theta+j))$ be independent random variables. A key fact is the representation

$$(19) \qquad K_{n,0} \stackrel{d}{=} (M_1 - 1)_+ + \cdots + (M_{n-1} - 1)_+ + M_n.$$

To prove this, one needs to set $M_j = \#\{k: S_k \in (E_{n-j,n}, E_{n-j+1,n})\}, j = 1, \ldots, n$ (with the convention $E_{0,n} = 0$), which is the number of points of a rate-$\theta$ Poisson process which fall between consecutive order statistics. The assertion about the joint distribution of $M_j$'s follows from the independence property of the Poisson process and the observation that the differences $E_{n,n} - E_{n-1,n}, E_{n-1,n} - E_{n-2,n}, \ldots, E_{n,1} - E_{n,0}$ are independent exponential variables with rates $1, 2, \ldots, n$. Now, counting the number of empty gaps $(S_k, S_{k+1})$ which fit in $(E_{n-j,n}, E_{n-j+1,n})$, we see that this is $M_n$ for $j = n$, and $(M_j - 1)_+$ for $j = 1, \ldots, n - 1$.

Passing to generating functions, (19) becomes

$$\mathbb{E}s^{K_{n,0}} = \frac{n}{n+\theta-\theta s} \prod_{j=1}^{n-1} \frac{j(j+2\theta-\theta s)}{(j+\theta)(j+\theta-\theta s)}$$

and (18) follows by sending $n \to \infty$ and evaluating the infinite product in terms of the gamma function [25]. The generating function of the stated mixed Poisson distribution is calculated by recalling that the generating function of Poisson($u$) is $e^{-u(1-s)}$ and that the Mellin transform of beta($1, \theta$) is $\mathbb{E}\xi^v = \theta B(\theta, 1 + v)$, whence

$$\mathbb{E}[\exp(\theta(1-s)\log \xi)] = \theta B(\theta, 1 + \theta - \theta s),$$

which is the same as (18). $\square$

The proof of Theorem 2.2 will exploit the poissonization technique, a well-known approach that goes back at least to Kac [22] (see also [12, 23] for the application of this technique to the balls-in-boxes scheme).

We shall first consider a sampling scheme in which exponential points $E_1, E_2, \ldots$ are thrown at the epochs of an independent Poisson process $\{\Pi(t) : t \geq 0\}$ with intensity one. After establishing convergence of $K_{\Pi(t), 0}$, we shall turn to that of $K_{n, 0}$.

PROOF OF THEOREM 2.2. (a) *Convergence in the Poisson model.* For $n, i \in \mathbb{N}_0$ and $t \geq 0$ set $a_n^{(i)} := \mathbb{P}\{K_{n, 0} = i\}$,

$$f^{(i)}(t) := \sum_{k=1}^{\infty} \frac{t^k}{k!} a_k^{(i)} \quad \text{and} \quad g^{(i)}(t) := e^{-t} f^{(i)}(t).$$

Notice that

$$g^{(0)}(t) + e^{-t} = \mathbb{P}\{K_{\Pi(t), 0} = 0\}, \qquad g^{(i)}(t) = \mathbb{P}\{K_{\Pi(t), 0} = i\}.$$

The equality of distributions (15) is equivalent to the following equalities:

$$a_0^{(0)} = 1, \qquad a_n^{(0)} = \sum_{k=0}^{n-1} a_k^{(0)} \mathbb{P}\{A_n^* = k\}, \qquad n \in \mathbb{N};$$

$$a_0^{(i)} = 0, \qquad a_n^{(i)} = a_n^{(i-1)} \mathbb{E}\bar{\xi}^n + \sum_{k=0}^{n-1} a_k^{(i)} \mathbb{P}\{A_n^* = k\}, \qquad i, n \in \mathbb{N},$$

from which we deduce after some calculations

$$g^{(0)}(t) = \mathbb{E}[g^{(0)}(t\bar{\xi})] + \mathbb{E}[e^{-t\bar{\xi}}] - e^{-t} - e^{-t}\mathbb{E}[f^{(0)}(t\bar{\xi})] =: \mathbb{E}[g^{(0)}(t\bar{\xi})] + f(t);$$

$$g^{(i)}(t) = \mathbb{E}[g^{(i)}(t\bar{\xi})] + e^{-t}(\mathbb{E}[f^{(i-1)}(t\bar{\xi})] - \mathbb{E}[f^{(i)}(t\bar{\xi})]), \qquad i \in \mathbb{N}.$$

Fix any $t_0 \in \mathbb{R}$ and define

$$f_1(t) := \mathbf{1}(t > t_0)(\mathbb{E}[\exp(-e^t \bar{\xi})] - \exp(-e^t)),$$
$$f_2(t) := \mathbf{1}(t \leq t_0)(\mathbb{E}[\exp(-e^t \bar{\xi})] - \exp(-e^t)),$$
$$f_3(t) := \mathbf{1}(t > t_0)\exp(-e^t)\mathbb{E}[f^{(0)}(e^t \bar{\xi})],$$
$$f_4(t) := \mathbf{1}(t \leq t_0)\exp(-e^t)\mathbb{E}[f^{(0)}(e^t \bar{\xi})].$$

Since $g^{(0)}$ is bounded and $g^{(0)}(0) = 0$,

$$g^{(0)}(e^t) = \int_{\mathbb{R}} f(e^{t-u}) \, d\left(\sum_{n=0}^{\infty} \mathbb{P}\{S_n \leq u\}\right).$$

If it were shown that $f_j$, $j = 1, 2, 3, 4$, was directly Riemann integrable (dRi) on $\mathbb{R}$, then since $f(e^t) = f_1(t) + f_2(t) - f_3(t) - f_4(t)$, we could apply the key renewal theorem to conclude that

(20)
$$\lim_{t \to \infty} \mathbb{P}\{K_{\Pi(t),0} = 0\} = \lim_{t \to \infty} g^{(0)}(e^t)$$
$$= \frac{1}{\mu} \int_0^{\infty} \frac{f(t)}{t} \, dt = 1 - \frac{1}{\mu} \sum_{j=1}^{\infty} \frac{\mathbb{E}\bar{\xi}^j}{j} \mathbb{P}\{K_{j,0} = 0\}.$$

We will only prove that $f_3$ and $f_4$ are dRi, the analysis of $f_1$ and $f_2$ being similar. Since $f_3$ and $f_4$ are continuous and positive on the sets $\{t \leq t_0\}$ and $\{t > t_0\}$, respectively, it suffices to find dRi majorants. We have

$$f_3(t) \leq \mathbf{1}(t > t_0)(\mathbb{E}[\exp(-e^t(1 - \bar{\xi}))] - \exp(-e^t))$$
$$\leq \mathbf{1}(t > t_0)\mathbb{E}[\exp(-e^t(1 - \bar{\xi}))] =: f_5(t),$$
$$f_4(t) \leq \mathbf{1}(t \leq t_0)(\mathbb{E}[\exp(-e^t(1 - \bar{\xi}))] - \exp(-e^t))$$
$$\leq \mathbf{1}(t \leq t_0)(1 - \exp(-e^t)) =: f_6(t).$$

The functions $f_5$ and $f_6$ are dRi, since they are bounded, monotone on the sets $\{t \leq t_0\}$ and $\{t > t_0\}$, respectively, and integrable. Integrability of $f_5$ follows from the condition $\nu < \infty$. This completes the proof of (20).

Arguing in the same manner as for the case $i = 0$, we conclude that for $i \in \mathbb{N}$

$$\lim_{t \to \infty} \mathbb{P}\{K_{\Pi(t),0} = i\} = \lim_{t \to \infty} g^{(i)}(e^t) = \frac{1}{\mu} \int_0^{\infty} \frac{e^{-t}(\mathbb{E}f^{(i-1)}(t\bar{\xi}) - \mathbb{E}f^{(i)}(t\bar{\xi}))}{t} \, dt$$
$$= \frac{1}{\mu} \sum_{j=1}^{\infty} \frac{\mathbb{E}\bar{\xi}^j}{j}(\mathbb{P}\{K_{j,0} = i - 1\} - \mathbb{P}\{K_{j,0} = i\}).$$

Assume now that $\mu = \infty$ and $\nu < \infty$. It suffices to prove that, as $t \to \infty$, $g^{(0)}(t) = e^{-t}\sum_{k=0}^{\infty}\frac{t^k}{k!}a_k^{(0)} \to 1$. Since $g^{(0)}(0) = 1$, $g^{(0)}$ is bounded and satisfies

$$g^{(0)}(t) = \mathbb{E}[g^{(0)}(t\bar{\xi})] - e^{-t}\mathbb{E}[f^{(0)}(t\bar{\xi})],$$

we conclude that

$$g^{(0)}(e^t) = 1 - \int_{\mathbb{R}} \exp(e^{t-u})\mathbb{E}[f^{(0)}(e^{t-u}\bar{\xi})] \, d\left(\sum_{n=0}^{\infty} \mathbb{P}\{S_n \leq u\}\right).$$

In the same way as in the first part of the proof we check that the key renewal theorem applies to yield

$$\lim_{t\to\infty} g^{(0)}(e^t) = 1 - \frac{1}{\mu}\int_0^\infty \frac{e^{-u}\mathbb{E}[f^{(0)}(u\bar\xi)]}{u} du = 0$$

(the last integral converges in view of the condition $\nu < \infty$). Thus, we have already proved that, as $t \to \infty$, $K_{\Pi(t),0} \xrightarrow{d} K_{\infty,0}$. Notice that, if $\mu < \infty$, then

$$\mathbb{E}K_{\infty,0} = \sum_{i=1}^\infty \mathbb{P}\{K_{\infty,0} \geq i\} = \frac{1}{\mu}\sum_{i=1}^\infty\sum_{j=1}^\infty \frac{\mathbb{E}\bar\xi^j}{j}\mathbb{P}\{K_{j,0} = i-1\} = \frac{1}{\mu}\sum_{j=1}^\infty \frac{\mathbb{E}\bar\xi^j}{j} = \frac{\nu}{\mu}.$$

(b) *Depoissonization.* For any fixed $\varepsilon \in (0,1)$ and $x > 0$, we have

$$\mathbb{P}\{K_{\Pi(t),0} > x\}$$
$$\leq \mathbb{P}\{K_{\Pi(t),0} > x, \lfloor(1-\varepsilon)t\rfloor \leq \Pi(t) \leq \lfloor(1+\varepsilon)t\rfloor\} + \mathbb{P}\{|\Pi(t) - t| > \varepsilon t\}$$
$$\leq \mathbb{P}\Big\{\max_{\lfloor(1-\varepsilon)t\rfloor \leq i \leq \lfloor(1+\varepsilon)t\rfloor} K_{i,0} > x\Big\} + \mathbb{P}\{|\Pi(t) - t| > \varepsilon t\}$$
$$= \mathbb{P}\{K_{\lfloor(1-\varepsilon)t\rfloor,0} > x\} + \mathbb{P}\Big\{K_{\lfloor(1-\varepsilon)t\rfloor,0} \leq x, \max_{\lfloor(1-\varepsilon)t\rfloor+1 \leq i \leq \lfloor(1+\varepsilon)t\rfloor} K_{i,0} > x\Big\}$$
$$+ \mathbb{P}\{|\Pi(t) - t| > \varepsilon t\} := I_1(t) + I_2(t) + I_3(t).$$

Similarly,

(21)
$$\mathbb{P}\{K_{\Pi(t),0} \leq x\}$$
$$\leq \mathbb{P}\{K_{\lfloor(1+\varepsilon)t\rfloor,0} \leq x\}$$
$$+ \mathbb{P}\Big\{K_{\lfloor(1+\varepsilon)t\rfloor,0} > x, \min_{\lfloor(1-\varepsilon)t\rfloor \leq i \leq \lfloor(1+\varepsilon)t\rfloor-1} K_{i,0} \leq x\Big\}$$
$$+ \mathbb{P}\{|\Pi(t) - t| > \varepsilon t\} := J_1(t) + J_2(t) + I_3(t).$$

If exponential points $E_{\lfloor(1-\varepsilon)t\rfloor+1}, \ldots, E_{\lfloor(1+\varepsilon)t\rfloor}$ fall to the left from the point $E_{\lfloor(1-\varepsilon)t\rfloor,\lfloor(1-\varepsilon)t\rfloor}$, then

$$\max_{\lfloor(1-\varepsilon)t\rfloor+1 \leq i \leq \lfloor(1+\varepsilon)t\rfloor} K_{i,0} \leq K_{\lfloor(1-\varepsilon)t\rfloor,0}$$

and also

$$K_{\lfloor(1+\varepsilon)t\rfloor,0} \leq \min_{\lfloor(1-\varepsilon)t\rfloor \leq i \leq \lfloor(1+\varepsilon)t\rfloor-1} K_{i,0},$$

which means that neither the event defining $I_2(t)$, nor $J_2(t)$, can hold. Therefore,

$$\max(I_2(t), J_2(t)) \leq \mathbb{P}\Big\{\max_{\lfloor(1-\varepsilon)t\rfloor+1 \leq i \leq \lfloor(1+\varepsilon)t\rfloor} E_i > E_{\lfloor(1-\varepsilon)t\rfloor,\lfloor(1-\varepsilon)t\rfloor}\Big\}$$

$$= \mathbb{E}(1 - (1 - e^{-E_{\lfloor(1-\varepsilon)t\rfloor,\lfloor(1-\varepsilon)t\rfloor}})^{\lfloor(1+\varepsilon)t\rfloor - \lfloor(1-\varepsilon)t\rfloor})$$

$$= 1 - \frac{\lfloor(1-\varepsilon)t\rfloor}{\lfloor(1+\varepsilon)t\rfloor}.$$

By a large deviation result (see, e.g., [2]), there exist positive constants $\delta_1$ and $\delta_2$ such that, for all $t > 0$,

$$I_3(t) \leq \delta_1 e^{-\delta_2 t}.$$

Select now $t$ such that $(1 - \varepsilon)t = n \in \mathbb{N}$. Then from the calculations above we get

$$\mathbb{P}\{K_{\Pi(n/(1-\varepsilon)),0} > x\}$$
$$\leq \mathbb{P}\{K_{n,0} > x\} + 1 - n/[(1+\varepsilon)n/(1-\varepsilon)] + \delta_1 \exp^{-\delta_2 n/(1-\varepsilon)}.$$

Sending first $n \uparrow \infty$ and then $\varepsilon \downarrow 0$, we obtain

$$\liminf_{n\to\infty} \mathbb{P}\{K_{n,0} > x\} \geq \mathbb{P}\{K_{\infty,0} > x\}$$

at all continuity points $x$ of the right-hand side. The same argument applied to (21) establishes the converse inequality for the upper limit.

(c) *Convergence of the mean value.* Denote by $H(x) := \sum_{k=0}^\infty \mathbb{P}\{S_k \leq x\}$ the renewal function and notice that

(22) $$\int_0^\infty e^{-sx} dH(x) = \frac{1}{1 - \mathbb{E}\bar\xi^s}, \qquad s > 0.$$

We have

$$\mathbb{E}K_{n,0} = \mathbb{E}\left[\sum_{k=0}^\infty ((1 - e^{-S_k} + e^{-S_{k+1}})^n - (1 - e^{-S_k})^n)\right]$$

$$= \int_0^\infty (\mathbb{E}(1 - e^{-x}\xi)^n - (1 - e^{-x})^n) dH(x)$$

$$= \int_0^\infty \left(\sum_{k=1}^n (-1)^{k+1} \binom{n}{k} e^{-kx}(1 - \mathbb{E}\xi^k)\right) dH(x)$$

$$\stackrel{(22)}{=} \sum_{k=1}^n (-1)^{k+1} \binom{n}{k} \frac{1 - \mathbb{E}\xi^k}{1 - \mathbb{E}\bar\xi^k}.$$

The conditions (5) imply that $\mu < \infty$ and $\nu < \infty$. The relation

(23) $$\lim_{n\to\infty} \mathbb{E}K_{n,0} = \nu/\mu$$

follows by an application of [9], Theorem 2(ii), to the formula for $\mathbb{E}K_{n,0}$. The cited result relies upon complex analysis and requires a sufficient large domain of definition of the Mellin transforms of $\xi$ and $\bar\xi$, which is secured by our assumption (5). We perceive that (23) holds whenever $\nu < \infty$, but have no proof of this conjecture so far.

For $n \in \mathbb{N}_0$ set $\varkappa_n := \mathbb{E} K_{n,0}$, $r_n := \frac{\mathbb{E} \bar{\xi}^n}{1 - \mathbb{E} \bar{\xi}^n}$. With decrement matrix as in (8), we have, according to (15),

$$\varkappa_0 = 0, \qquad \varkappa_n = \sum_{m=1}^{n} q(n:m) \varkappa_{n-m} + r_n, \qquad n \in \mathbb{N},$$

which is of the same form as [10], (11). Then $\varkappa_n$ is given by

$$\varkappa_n = \sum_{m=1}^{n-1} g(n,m) r_m + r_n, \qquad n \in \mathbb{N}$$

(compare to [10], (12)), where $g(n,m)$ was defined on page 10. Assuming that $\nu = \infty$ and $\mu < \infty$ and using (16) along with Fatou's lemma gives

$$\liminf_{n \to \infty} \varkappa_n \geq \sum_{m=1}^{\infty} \frac{1 - \mathbb{E} \bar{\xi}^m}{\mu m} \frac{\mathbb{E} \bar{\xi}^m}{1 - \mathbb{E} \bar{\xi}^m} = \frac{1}{\mu} \sum_{m=1}^{\infty} \frac{\mathbb{E} \bar{\xi}^m}{m} = \infty,$$

where the last series diverges in view of the condition $\nu = \infty$. $\square$

Exactly the same argument as above can be exploited for proving that $K_{n,0} + \cdots + K_{n,r}$ converges in distribution. However, for $r \geq 2$ calculations get complicated and in Proposition 5.2 we content ourselves with the case $r = 1$.

PROPOSITION 5.2. *If $\nu < \infty$, then, as $n \to \infty$, $K_{n,0} + K_{n,1}$ converges in distribution to a random variable $K_{01}$. If also $\mu < \infty$, then*

$$\mathbb{P}\{K_{01} \geq 1\} = \frac{1}{\mu} \left( \mathbb{E} \xi + \sum_{j=2}^{\infty} \left( \frac{\mathbb{E} \bar{\xi}^j}{j} + \mathbb{E} \bar{\xi}^j \xi \right) \mathbb{P}\{K_{j,0} + K_{j,1} = 0\} \right),$$

$$\mathbb{P}\{K_{01} \geq i\} = \frac{1}{\mu} \Bigg( \mathbb{E} \bar{\xi}^{i-2} \mathbb{E} \xi (1 - \mathbb{E} \xi^2)$$

$$+ \sum_{j=2}^{\infty} \left( \frac{\mathbb{E} \bar{\xi}^j}{j} + \mathbb{E} \bar{\xi}^j \xi \right) \mathbb{P}\{K_{j,0} + K_{j,1} = i - 1\} \Bigg),$$

$$i = 2, 3, \ldots,$$

*and $\mathbb{E} K_{01} = (\nu + 1)/\mu$, but if $\mu = \infty$, then $K_{01} = 0$ a.s.*

SKETCH OF THE PROOF. For $n, i \in \mathbb{N}_0$ and $t \geq 0$ set $Y_n := K_{n,0} + K_{n,1}$,

$$f^{(i)}(t) = \sum_{k=2}^{\infty} \frac{t^k}{k!} \mathbb{P}\{Y_k = i\} \quad \text{and} \quad g^{(i)}(t) := e^{-t} f^{(i)}(t).$$

Use the recursion

$$Y_0 = 0, \qquad Y_n \stackrel{d}{=} Y_{A_n^*} + \mathbf{1}(A_n^* \geq n-1), \qquad n \in \mathbb{N}, \tag{24}$$

where $A_n^*$ is independent of $\{Y_k : k \in \mathbb{N}\}$, to obtain

$$\begin{aligned}
g^{(0)}(t) &= \mathbb{E}g^{(0)}(t\bar\xi) + \mathbb{E}e^{-t\bar\xi} - e^{-t} - e^{-t}\mathbb{E}f^{(0)}(t\bar\xi) \\
&\quad - te^{-t}\mathbb{E}\xi(f^{(0)}(t\bar\xi) + 1);
\end{aligned} \tag{25}$$

$$\begin{aligned}
g^{(i)}(t) &= e^{-t}\sum_{n=2}^{\infty}\frac{t^n \mathbb{P}\{Y_{n-1} = i-1\}\mathbb{E}\bar\xi^{n-1}\xi}{(n-1)!} + e^{-t}g^{(i-1)}(t) \\
&\quad + \mathbb{E}(f^{(i)}(t\bar\xi) + t\bar\xi(\mathbb{E}\bar\xi)^{i-1}\mathbb{E}\xi)(e^{-t\bar\xi} - e^{-t} - te^{-t}\xi).
\end{aligned} \tag{26}$$

In the same way as in the proof of Theorem 2.2, we can justify using the key renewal theorem in (25) and (26) to get a poissonized version of the result. Our depoissonization argument used in the proof of Theorem 2.2 works without changes. Justification of the only step that may require explanation is as follows: the inequality $K_{n,0} + K_{n,1} < K_{m,0} + K_{m,1}$, $n < m$ implies that at least one of the exponential points $E_{n+1}, \ldots, E_m$ falls to the right from $E_{n,n}$. $\square$

**6. Proof of Theorem 2.3.** Assume that $\nu < \infty$ and that $(K_n^* - b_n)/a_n$ converges in distribution to a random variable $X$ with some proper and nondegenerate probability law. According to Theorem 2.1, the latter can occur if one of five conditions (1)–(5) holds and also $a_n \to \infty$. Notice that for each of these conditions there exist distributions that satisfy it together with the condition $\nu < \infty$ [as, e.g., in Examples 2.5 and 2.6]. By Theorem 2.2(a) and the Markov inequality, $K_{n,0}/a_n = (K_n^* - K_n)/a_n$ goes to 0 in probability. Therefore, $(K_n - b_n)/a_n \stackrel{d}{\to} X$. Using Proposition 5.2 and the Markov inequality, we conclude that $(K_n - K_{n,1} - b_n)/a_n \stackrel{d}{\to} X$.

To prove the result for $W_n$, consider (13) and assume that there exists an increasing sequence $\{n_k : k \in \mathbb{N}\}$ such that $V_{n_k} \to \infty$. According to Theorem 2.2, we get from (13) $K_{\infty,0} \stackrel{d}{=} K_{\infty,0} + 1$, which is absurd. Thus, the sequence $\{V_n : n \in \mathbb{N}\}$ is bounded in probability, which implies that

$$\{U_{V_n} - \mathbf{1}(V_n = 0) : n \in \mathbb{N}\} \text{ is bounded in probability,} \tag{27}$$

as well. An appeal to (14) allows us to conclude that $(W_n - b_n)/a_n \stackrel{d}{\to} X$.

Assume now that either $(K_n - b_n)/a_n$, $(K_n - K_{n,1} - b_n)/a_n$, or $(W_n - b_n)/a_n$ converges in distribution to a random variable $X$ with some proper and nondegenerate probability distribution. Essentially in the same way as for $K_n^*$ [but now using also either the result of Theorem 2.2 or Proposition 5.2 or (27)], we can prove that $a_n \to \infty$, and the same argument as above proves that $(K_n^* - b_n)/a_n \stackrel{d}{\to} X$. The proof is complete.

## 7. Proof of Theorem 2.4.

*Case $\mu < \infty$.* With $g(n,m)$ defined on page 10, we have

$$\mathbb{P}\{Z_n = m\} = g(n,m)\mathbb{P}\{A_m = 0\}.$$

Since, according to (8),

$$\mathbb{P}\{A_m = 0\} = \frac{\mathbb{E}\xi^m}{1 - \mathbb{E}\bar{\xi}^m},$$

an appeal to (16) completes the proof of this case.

*Case $\mu = \infty$.* Denote $\widetilde{U}(z) := \inf\{z - S_n : S_n < z, n \in \mathbb{N}_0\}$ the undershoot at $z > 0$. For $k \in \{1, 2, \ldots, n\}$ we have

(28) $$\mathbb{P}\{Z_n > k\} = \mathbb{P}\{\widetilde{U}(E_{n,n}) > E_{n,n} - E_{n-k,n}\}.$$

Assume first that $\alpha \in [0, 1)$ and for fixed $\varepsilon \in (0, 1)$ set $k_n := \lfloor n^\varepsilon \rfloor$. Since $E_{n,n}$ is independent of the undershoot and tends to $+\infty$ in probability, an appeal to [6], Theorem 8.6.3, allows us to conclude that

$$\frac{\widetilde{U}(E_{n,n})}{E_{n,n}} \xrightarrow{d} \widehat{Z}_\alpha,$$

where the distribution of $\widehat{Z}_0$ is $\delta_1$, degenerate at point 1, and for $\alpha \in (0, 1)$, $\widehat{Z}_\alpha$ has the beta $(1-\alpha, \alpha)$ distribution. Using the convergence of $E_{n.n} - \log n$, we obtain $E_{n,n}/\log n \to 1$ in probability. Since, for $x > 0$,

$$\mathbb{P}\{E_{n,n} - E_{n-k_n,n} \leq x\} = (1 - e^{-x})^{k_n},$$

we can easily check that $(E_{n,n} - E_{n-k_n,n})/\log n \xrightarrow{d} \varepsilon$. Therefore,

$$\frac{E_{n,n} - E_{n-k_n,n}}{E_{n,n}} \xrightarrow{d} \varepsilon.$$

Now the result follows from the relation

$$\mathbb{P}\left\{\frac{\log Z_n}{\log n} > \varepsilon\right\} = \mathbb{P}\{Z_n > k_n\} \stackrel{(28)}{=} \mathbb{P}\left\{\frac{\widetilde{U}(E_{n,n})}{E_{n,n}} > \frac{E_{n,n} - E_{n-k_n,n}}{E_{n,n}}\right\}$$

$$\to \mathbb{P}\{\widehat{Z}_\alpha > \varepsilon\}.$$

Indeed, while in the case $\alpha \in (0, 1)$, each $\varepsilon \in (0, 1)$ is a continuity point of the distribution of $\widehat{Z}_\alpha$, in the case $\alpha = 0$ the relation establishes the convergence in probability $\log Z_n / \log n \to 1$ (notice that $\log Z_n / \log n \leq 1$ a.s.).

Consider now the remaining case $\alpha = 1$. For fixed $\varepsilon \in (0, 1)$ set $k_n := \lfloor \exp(m^{-1}(\varepsilon m(\log n))) \rfloor$, where $m^{-1}(\cdot)$ is the increasing and continuous inverse of $m(x) = \int_0^x \mathbb{P}\{-\log \bar{\xi} > y\}\, dy$, $x > 0$. Using again the independence of $E_{n,n}$ and the undershoot and exploiting [7], Theorem 6, leads to the conclusion

$$\frac{m(\widetilde{U}(E_{n,n}))}{m(E_{n,n})} \xrightarrow{d} \widetilde{Z},$$

where the law of $\widetilde{Z}$ is uniform$[0,1]$. Fix any $i \in \mathbb{N}$. It is well known that $m(x)$ is slowly varying at $\infty$. Therefore, $m^i(\log x)$ is also slowly varying at $\infty$, and as $s \downarrow 0$, $m^i(-\log(1-e^{-s})) \sim m^i(-\log s)$, where $f \sim g$ means that the ratio $f/g$ goes to one. Applying Proposition 1.5.8 and Theorem 1.7.1' from [6] to the equality

$$\mathbb{E}[m^i(E_{n,n})] = n\int_0^\infty m^i(-\log(1-e^{-s}))e^{-ns}\,ds$$

we get $\mathbb{E}[m^i(E_{n,n})] \sim m^i(\log n)$. Similarly,

$$\mathbb{E}m^i(E_{n,n} - E_{n-k_n,n}) \sim m^i(\log k_n)$$
$$\sim m^i(\log \exp(m^{-1}(\varepsilon m(\log n)))) = \varepsilon^i m^i(\log n).$$

The last two relations (with $i=1$ and $i=2$) together with Chebyshev's inequality imply that

$$m(E_{n,n})/m(\log n) \xrightarrow{d} 1 \quad \text{and} \quad m(E_{n,n} - E_{n-k_n,n})/m(\log n) \xrightarrow{d} \varepsilon.$$

Consequently, $m(E_{n,n} - E_{n-k_n,n})/m(E_{n,n}) \xrightarrow{d} \varepsilon$. To finish the proof, it remains to note that

$$\mathbb{P}\left\{\frac{m(\log Z_n)}{m(\log n)} > \varepsilon\right\} = \mathbb{P}\{Z_n > k_n\}$$
$$\stackrel{(28)}{=} \mathbb{P}\left\{\frac{m(\widetilde{U}(E_{n,n}))}{m(E_{n,n})} > \frac{m(E_{n,n} - E_{n-k_n,n})}{m(E_{n,n})}\right\}$$
$$\to \mathbb{P}\{\widetilde{Z} > \varepsilon\}.$$

## APPENDIX

Asymptotic behavior of the first passage time processes $N_t = \inf\{k \geq 1 : S_k \geq t\}$, for $S_k = X_1 + \cdots + X_k$ a zero-delayed random walk with positive steps, was investigated by many authors (see, e.g., [5, 8, 18, 19]). The next proposition is a summary of results scattered in the literature.

PROPOSITION A.1. *Assume that $X_1 > 0$ a.s. and that the distribution of $X_1$ is nonlattice. The following assertions are equivalent:*

(i) *There exist functions $a(t) > 0, b(t) \in \mathbb{R}$ such that, as $t \to \infty$, $(N_t - b(t))/a(t)$ converges weakly to a nondegenerate and proper probability law.*

(ii) *Either the distribution of $X_1$ belongs to the domain of attraction of a stable law, or $\mathbb{P}\{X_1 > x\}$ slowly varies at $\infty$.*

Set $\mu = \mathbb{E}X_1$ and $\sigma^2 = \mathbb{D}X_1$.

(a) *If $\sigma^2 < \infty$, then, with $b(t) = \mu^{-1}t$ and $a(t) = (\mu^{-3}\sigma^2 t)^{1/2}$, the limiting law is standard normal.*

(b) *If $\sigma^2 = \infty$ and*
$$\int_0^x y^2 \mathbb{P}\{X_1 \in dy\} \sim L(x) \qquad \text{as } x \to \infty,$$
*for some $L$ slowly varying at $\infty$, then, with $b(t) = \mu^{-1}t$ and $a(t) = \mu^{-3/2}c(t)$, where $c(t)$ is any function satisfying*
$$\lim_{t \to \infty} tL(c(t))/c^2(t) = 1,$$
*the limiting law is standard normal.*

(c) *Assume that the relation*
$$\mathbb{P}\{X_1 > x\} \sim x^{-\alpha} L(x) \qquad \text{as } x \to \infty, \tag{29}$$
*where $L$ is some function slowly varying at $\infty$, holds with $\alpha \in [1, 2)$, and that in the case $\alpha = 1$ also $\mu < \infty$. Then, with $b(t) = \mu^{-1}t$ and $a(t) = \mu^{-(\alpha+1)/\alpha}c(t)$, where $c(t)$ is any function satisfying*
$$\lim_{t \to \infty} tL(c(t))/c^\alpha(t) = 1,$$
*the limiting law is $\alpha$-stable with characteristic function*
$$t \mapsto \exp\{-|t|^\alpha \Gamma(1-\alpha)(\cos(\pi\alpha/2) + i\sin(\pi\alpha/2)\operatorname{sgn}(t))\}, \qquad t \in \mathbb{R}.$$

(d) *Assume that $\mu = \infty$ and the relation (29) holds with $\alpha = 1$. Let $c$ be any positive function satisfying $\lim_{x \to \infty} xL(c(x))/c(x) = 1$ and set $\psi(x) := x \int_{\exp(-c(x))}^1 \mathbb{P}\{X_1 \leq y\} y^{-1} dy$. Let $b(t)$ be any positive function satisfying $b(\psi(t)) \sim \psi(b(t)) \sim t$. Then, with $a(t) = b(t)c(b(t))/t$, the limiting law is 1-stable with characteristic function*
$$t \mapsto \exp\{-|t|(\pi/2 - i\log|t|\operatorname{sgn}(t))\}, \qquad t \in \mathbb{R}.$$

(e) *If the relation (29) holds with $\alpha \in [0, 1)$, then, with $b(t) \equiv 0$ and $a(t) = t^\alpha/L(t)$, the limiting law $\theta_\alpha$ is a scaled Mittag–Leffler (exponential, if $\alpha = 0$) with moments*
$$\int_0^\infty x^n \theta_\alpha(dx) = \frac{n!}{\Gamma^n(1-\alpha)\Gamma(1+n\alpha)}, \qquad n \in \mathbb{N}.$$

All the above statements remain valid if the continuous variable $t$ is replaced by discrete variable $\log n$, $n \in \mathbb{N}$, as has been used in this paper.

In the lattice case, when $X_1$ assumes only positive integer values, the whole range of possible distributional limits follow from Theorems 1.2, 1.5 and Proposition 3.1 in [21]. Although the first two of these results were formulated for other variables, they apply to $N_t$ as well. As in [21], the same asymptotic results are readily extendible to the first passage time processes for random walks with positive nonlattice steps.

A. V. GNEDIN
MATHEMATICAL INSTITUTE
UTRECHT UNIVERSITY
POSTBUS 80010
3508 TA UTRECHT
THE NETHERLANDS
E-MAIL: A.V.Gnedin@math.uu.nl

A. M. IKSANOV
P. NEGADAJLOV
FACULTY OF CYBERNETICS
NATIONAL T. SHEVCHENKO UNIVERSITY
01033 KIEV
UKRAINE
E-MAIL: iksan@unicyb.kiev.ua
npasha@ukr.net

U. RÖSLER
MATHEMATISCHES SEMINAR
CHRISTIAN–ALBRECHTS UNIVERSITÄT ZU KIEL
LUDEWIG–MEYN STR. 4
D-24098 KIEL
GERMANY
E-MAIL: roesler@computerlabor.math.uni-kiel.de